\documentclass{tac}
\usepackage{amssymb}
\usepackage{amsmath}
\usepackage[cmtip,all]{xypic}

\author{Daniel G. Davis} 
\address{Department of Mathematics\\
Wesleyan University\\
265 Church St.\\
Middletown, CT 06459-0128}

\title{The site $R^+_G$ for a profinite group $G$}
\copyrightyear{2006}
\eaddress{dgdavis@wesleyan.edu}
 
\let\pf\proof
\let\epf\endproof

\newcommand{\G}{G}
\newcommand{\R}{R^+_{\G}}
\newcommand{\spc}{\!}
\DeclareMathOperator*{\holim}{holim}

\def\:{\colon}

\begin{document}
\maketitle
\begin{abstract}
Let $G$ be a non-finite profinite group and let 
$G-\mathbf{Sets}_{df}$ be the canonical site of finite 
discrete $G$-sets. Then the category $R^+_G$, defined by 
Devinatz and Hopkins, is the category obtained by 
considering $G-\mathbf{Sets}_{df}$ together with the 
profinite $G$-space $G$ itself, with morphisms being 
continuous $G$-equivariant maps. We show that $R^+_G$ 
is a site when equipped with the pretopology of epimorphic 
covers. Also, we explain why the associated topology on $R^+_G$ 
is not 
subcanonical, and hence, not canonical. We note that, since 
$R^+_G$ is a site, there is 
automatically a model category structure on the category of 
presheaves of spectra on the site. Finally, we point out that such 
presheaves of spectra are a nice way of organizing the data that 
is obtained by taking the homotopy fixed points of a continuous 
$G$-spectrum with respect to the open subgroups of $G$. 
\end{abstract} 
\section{Introduction}
Let $G$ be a profinite group that is not a finite group. Let $\R$ 
be the category with objects all finite
discrete left $G$-sets together with the left $\G$-space $G$. The
morphisms of $R^+_G$ are the continuous $G$-equivariant maps. 
Since $G$ is not 
finite, the object $G$ in $\R$ is very different in character from all the other 
objects of $\R$. In this paper, we 
show that $\R$ is a site when equipped with the pretopology of epimorphic covers. 
\par
As far as the author knows, the category $R^+_G$ is first defined and used in the paper \cite{DH}, by Ethan Devinatz and Mike Hopkins. Let $G_n$ be the profinite group $S_n \rtimes \mathrm{Gal}(\mathbb{F}_{p^n}/\mathbb{F}_p),$ where $S_n$ is the $n$th Morava stabilizer group. In \cite[Theorem 1]{DH}, Devinatz and Hopkins construct a contravariant functor - that is, a presheaf - 
\[\mathbf{F} \: (R^+_{G_n})^\mathrm{op} \rightarrow (\mathcal{E}_{\infty})_{K(n)},\] to the category $(\mathcal{E}_{\infty})_{K(n)}$ of $K(n)$-local commutative $S$-algebras (see \cite{EKMM}), where $K(n)$ is the $n$th Morava $K$-theory (see \cite[Chapter 9]{Rudyak} for an exposition of $K(n)$). The functor $\mathbf{F}$ has the properties that, if $U$ is an open subgroup of $G_n$, then $\mathbf{F}(G_n/U) = E_n^{dhU}$, and $\mathbf{F}(G_n)=E_n$, where $E_n$ is the $n$th Lubin-Tate spectrum (for salient facts about $E_n$ and its importance in homotopy theory, see 
\cite[Introduction]{DevinatzHopkins}), and $E_n^{dhU}$ is a spectrum that behaves like the $U$-homotopy fixed point spectrum of $E_n$ with respect to the continuous $U$-action. Since 
$\mathrm{Hom}_{R^+_{G_n}}(G_n,G_n) \cong G_n$, functoriality implies that $G_n$ acts 
on $E_n$ by maps of commutative $S$-algebras. In Section \ref{presheaves}, we will give several related examples of presheaves of spectra that illustrate the utility of the category $\R$. 
\par
The pretopology of epimorphic covers on a small category
$\mathcal{C}$ is the pretopology $\mathcal{K}$ given by all covering families
$\{f_i \: C_i \rightarrow C |\, i \in I\}$ such that $\phi \:
\coprod_{i \in I} C_i \rightarrow C$ is onto, where $C_i, C \in
\mathcal{C}$, $f_i \in \mathrm{Mor}_{\mathcal{C}}(C_i, C)$, and $I$ is some
indexing set. (Of course, one must prove that these covering families
actually give a pretopology on $\mathcal{C}$.) We note that we 
do not require that $\phi$ be a
morphism in $\mathcal{C}$; for our purposes, $\mathcal{C} = \R$ and we
only require that $\phi$ be an epimorphism in the category of all
$\G$-sets (so that $\phi$ does not have to be continuous). This assumption is
important for our work, since, for example, $G \coprod G$ is not in
$\R$. 
\par
The pretopology $\mathcal{K}$ is a familiar one. For example, for a profinite group
$G$, $\mathcal{K}$ is the standard basis used for the site $G-\mathbf{Sets}_{df}$ of
finite discrete $G$-sets (\cite[pg. 206]{Jardine}). However, there is an
important difference between $\R$ and $G-\mathbf{Sets}_{df}$: the latter category is closed under pullbacks, but it is easy to see that $\R$
does not have all pullbacks (this point will be discussed
later).  But in a category with pullbacks, the canonical
topology, the finest topology in which every representable presheaf is
a sheaf, is given by all covering families of universal effective
epimorphisms (see Expose IV, 4.3 of \cite{Demazure}). 
This implies that $G-\mathbf{Sets}_{df}$ is a site with the canonical topology when equipped with pretopology
$\mathcal{K}$. However, due to the lack of sufficient pullbacks, we cannot 
conclude that $\mathcal{K}$ gives $\R$ the
canonical topology. In fact, we will show that $\mathcal{K}$ does not generate the canonical topology, since $\mathcal{K}$ does not yield a subcanonical topology.
\par
Note that $R^+_{G}$ is built out of the two subcategories
$G \negthinspace - \negthinspace \mathbf{Sets}_{df}$ and the groupoid $G$. Since each of these categories is
a site via $\mathcal{K}$ (for $G$, this is verified in Lemma 2 below), it is natural to think 
that $\R$ is also a site via
$\mathcal{K}$. Our main result (Theorem \ref{site1}), verifies that this is indeed the case.
\par
As discussed earlier, $\mathbf{F}$ is a presheaf of spectra on the site $R^+_G$. More generally, there is the category $\mathrm{PreSpt}(R^+_G)$ of presheaves of spectra on $R^+_G$. Furthermore, since $R^+_G$ is a site, the work of Jardine (e.g., \cite{Jardinecanada}, \cite{Jardine}) implies that $\mathrm{PreSpt}(R^+_G)$ is a model category. We recall the definition of this model category in Section 5.
\par 
In \cite{cts}, the author showed that, given a continuous $G$-spectrum $Z$, then, for any open 
subgroup $U$ of $G$, there is a homotopy fixed point spectrum $Z^{hU}$, defined with respect to the continuous action of $U$ on $Z$. In Examples \ref{exone} and \ref{extwo}, we see that there is a presheaf that organizes in a functorial way the following data: $Z$, $Z^{hU}$ for all $U$ open in $G$, and the maps between these spectra that are induced by continuous $G$-equivariant maps between the $G$-spaces $G$ and $G/U$. Thus, $\mathrm{PreSpt}(R^+_G)$ is a natural category within which to work with continuous $G$-spectra. It is our hope that the model category structure on $\mathrm{PreSpt}(R^+_G)$ can be useful for the theory of homotopy fixed points for profinite groups, though we have not yet found any such applications.
\par
\noindent
\textbf{Acknowledgements.} When I first tried to make $\R$ a site, and was focusing on an abstract way of doing this, Todd Trimble helped me get started by
suggesting that I extend $\mathcal{K}$
to all of $\R$. Also, I thank him for pointing out Lemma \ref{one}. I thank Paul Goerss 
for discussions about this material. Also, I appreciate various conversations with Christian Haesemeyer about this work.
\section{Preliminaries}
\par
Before we prove our main results, we first collect some easy facts
which will be helpful later. As stated in the Introduction, $G$ always refers to an infinite
profinite group. (If the profinite group $G$ is finite, then $\R = G \negthinspace - 
\negthinspace \mathbf{Sets}_{df}$
and there is nothing to prove.)
\lemma\label{one}
Let $f \: C \rightarrow \G$ be any morphism in $\R$ with $C \neq
\varnothing$.  Then $C = \G$.
\endlemma
\pf
Choose any $c \in C$ and let $f(c) = \gamma$. Choose any $\delta \in
\G$. Then \[\delta = (\delta \gamma^{-1})\gamma = (\delta \gamma^{-1})
\cdot f(c) = f((\delta \gamma^{-1}) \cdot c),\] by the
$\G$-equivariance of $f$. Thus, $f$ is onto and $|\mathrm{im}(f)| = \infty$, so
that $C$ cannot be a finite set.
\epf
\lemma 
For a topological group $G$, let $G$ be the groupoid with the single
object $G$ and morphisms the $G$-equivariant maps $G \rightarrow
G$ given by right multiplication by some element of $G$. Then $G$ is a site
with the pretopology $\mathcal{K}$ of epimorphic covers.
\endlemma
\pf
Any diagram $G \overset{f}{\rightarrow} G \overset{g}{\leftarrow} G$, where
$f$ and $g$ are given by multiplication by $\gamma$ and $\delta$,
respectively, can be completed to a commutative square \[\xymatrix{G
\rto^{f'} \dto^{g'} & G \dto^{g} \\ G \rto^{f} & G,\\}\] where $f'$
and $g'$ are given by multiplication by $\delta^{-1}$ and
$\gamma^{-1}$, respectively. This property suffices to show that $G$
is a site with the atomic topology, in which every sieve is a covering
sieve if and only if it is nonempty. It is easy to see that the only nonempty
sieve of $G$ is $\mathrm{Mor}_G(G,G)$ itself. Thus, the only covering sieve of
$G$ is the maximal sieve. Since every morphism of $G$ is a
homeomorphism, in the pretopology $\mathcal{K}$, the collection of
covers is exactly the collection of all
nonempty subsets of $\mathrm{Mor}_G(G,G)$. Then it is easy to check that $\mathcal{K}$ is the maximal
basis that generates the atomic topology.
\epf
Observe that if $f \: \G \rightarrow \G$ is a morphism in
$\R$, then by $\G$-equivariance, $f$ is the map given by
multiplication by $f(1)$ on the right. As mentioned earlier, we have 
\lemma 
The category $G \negthinspace - \negthinspace \mathbf{Sets}_{df}$, a full subcategory of $\R$, is closed
under pullbacks.
\endlemma
\pf
The pullback of a diagram in $G \negthinspace - \negthinspace \mathbf{Sets}_{df}$ is formed simply by regarding the diagram as being in the category $T_G$ of discrete $G$-sets. The category $T_G$ is closed under pullbacks, as explained in \cite[pg. 31]{MacLane}.
\epf
We recall the following useful result and its proof.
\lemma\label{identify}\label{four}
Let $X$ be any finite set in $\R$. We write $X = \coprod_{i=1}^n
\overline{x_i}$, the disjoint union of all the distinct orbits
$\overline{x_i}$, with each $x_i$ a representative. Then $X$ is
homeomorphic to $\coprod_{i=1}^n \G/U_i$, where $U_i={\G}_{x_i}$ is
the stabilizer in $\G$ of $x_i$.
\endlemma

\pf
Let $f \: \G/U_i \rightarrow \overline{x_i}$ be given by $f(\gamma
U_i) = \gamma \cdot x_i$. Since $X$ is a discrete $\G$-set, the
stabilizer $U_i$ is an open subgroup of $\G$ with finite index, so that $\G/U_i$
is a finite set. Then $f$ is open and continuous since it is a map between
discrete spaces. Also, it is clear that $f$ is onto. Now suppose $\gamma
U_i = \delta U_i$. Then $\gamma^{-1}\delta \in U_i$, so that
$(\gamma^{-1}\delta)\cdot x_i = (\gamma^{-1}) \cdot (\delta \cdot x_i)
= x_i$. Thus, $\gamma \cdot x_i = \delta \cdot x_i$ and $f$ is
well-defined. Assume that $\gamma \cdot x_i = \delta \cdot x_i$. Then
$\gamma^{-1}\delta \in {\G}_{x_i}$ so that $f$ is a monomorphism.
\epf

\lemma\label{five}
Let $X$ be a finite discrete $G$-set in $\R$ and let $\psi \: \G \rightarrow X$
be any $\G$-equivariant function. Then $\psi$ is a morphism in $\R$.
\endlemma
 
\pf
As in Lemma \ref{identify}, we identify $X$ with $\coprod_{i=1}^n \G/U_i$. Since $\psi$ is $\G$-equivariant and
$\psi(\gamma) = \gamma \cdot \psi(1)$, $\psi$ is determined by
$\psi(1)$. Let $\psi(1) = \delta U_j$ for some $\delta \in \G$ and some $j$. Then for any $\gamma$ in $\G$, \[\gamma U_j = (\gamma
\delta^{-1} \delta)U_j = (\gamma \delta^{-1}) \cdot \psi(1) =
\psi(\gamma \delta^{-1}),\] so that $\mathrm{im} \, \psi = \G/U_j.$ Since $X$
is discrete, $\psi$ is continuous, if, for any $x \in X$,
$\psi^{-1}(x)$ is open in $\G$. It suffices, by the identification, to
let $x=\gamma U_j$, for any $\gamma \in \G$. Then 
\begin{align*}
\psi^{-1}(\gamma
U_j) & = \{ \zeta \in \G | \, \psi(\zeta) = \gamma U_j \} = \{ \zeta
\in G | \, \zeta \cdot (\delta U_j) = \gamma U_j \} \\ & = \{ \zeta \in G | \, \delta^{-1}
\zeta^{-1} \gamma \in U_j \} = \gamma U_j \delta^{-1}.
\end{align*} 
Since $U_j$
is open and multiplication on the left or the right is always a
homeomorphism in a topological group, we see that $\psi^{-1}(x)$ is an
open set in $\G$. 
\epf

\section{The proof of the main theorem}
With these lemmas in hand, we are ready for

\theorem\label{site1} 
For any profinite group $G$, the category $R^+_{G}$ equipped with the pretopology $\mathcal{K}$ of
epimorphic covers is a small site.
\endtheorem

Before proving the theorem, we first make some remarks about pullbacks
in $\R$ and how this affects our proof. In a category $\mathcal{C}$
with sufficient pullbacks, to prove that a pretopology is given by a function $K$, which assigns
to each object $C$ a collection $K(C)$ of families of morphisms with
codomain $C$, one must prove the stability axiom, which says the
following: if $\{f_i \: C_i \rightarrow C | \, i \in I\} \in
K(C)$, then for any morphism $g \: D \rightarrow C$, the family of
pullbacks \[\{\pi_L \: D \times_C C_i \rightarrow D | \, i \in I\}
\in K(D).\]

Let us examine what this axiom would require of $\R$. 

\begin{example}
The map $G \rightarrow *$ forms a covering family
and so the stability axiom requires that $G {\times _ { \{ * \} } } G
= G \times G$ be in $\R$.
\end{example}
\begin{example}
Let $C$ be any finite discrete $G$-set with more than
one element and with trivial $G$-action, $g \: \G
\rightarrow C$ any morphism, and consider the cover \[\{f_i \: C_i
\rightarrow C | \, i \in I\} \in K(C),\] where $C_j = C$ and $f_j \:
C \rightarrow C$ is the morphism mapping $C$ to $g(1)$, for some
$j \in I$. Because the action is trivial, $f_j$ is
$G$-equivariant. There certainly exist covers of $C$ of this form, since one
could let $f_k = \mathrm{id}_C$, for some $k \neq j$ in $I$, and then let the
other $f_i$ be any morphisms with codomain $C$. Then the stability
axiom requires that $G \times_C C$ exists in $\R$, but this is impossible,
since \[G \times_C C = \{(\gamma, c) | \, g(\gamma) = f_j(c) \} =
\{(\gamma, c) | \, \gamma \cdot g(1) = g(1) \} = G_{g(1)} \times C = G
\times C.\]
\end{example}
\par
Thus, the stability axiom for a pretopology must be altered so that one
still obtains a topology. We list the correct axioms for our situation
below. They are taken from \cite[Exercise 3, pg. 156]{MacLane}.
 
\begin{enumerate}
\item
If $f \: C' \rightarrow C$ is an isomorphism, then $\{f \: C'
\rightarrow C\} \in K(C)$.
\item
(stability axiom) If $\{f_i \: C_i \rightarrow C |\, i \in I\} \in
K(C)$, then for any morphism $g \: D \rightarrow C$, there exists
a cover $\{h_j \: D_j \rightarrow D |\, j \in J\} \in K(D)$ such
that for each $j$, $g \circ h_j$ factors through some $f_i$.
\item
(transitivity axiom)
If $\{f_i \: C_i \rightarrow C |\, i \in I\} \in K(C)$, and if for
each $i \in I$ there is a family $\{g_{ij} \: D_{ij} \rightarrow
C_i |\, j\in I_i \} \in K(C_i)$, then the family of composites \[\{f_i
\circ g_{ij} \: D_{ij} \rightarrow C |\, i\in I, j\in I_i \}\] is
in $K(C)$. 
\end{enumerate}

\pf[of Theorem \ref{site1}]
It is clear that the pretopology of epimorphic covers satisfies
axiom (1) above. Also, it is easy to see that axiom (3) holds. Indeed, using the
above notation, choose any $c \in C$. Then there is some $c_i \in C_i$
for some $i$, such that $f_i(c_i)= c$. Similarly, there must be some
$d_{ij} \in D_{ij}$ for some $j$, such that $g_{ij}(d_{ij}) =
c_i$. Hence, $(f_i \circ g_{ij})(d_{ij}) = f_i(c_i) = c$, so that
$\coprod_{i,j} D_{ij} \rightarrow C$ is onto. This verifies (3).
We verify (2) by considering five cases.

{\it Case $\mathrm{(}$1$\, \mathrm{)}$}: Suppose that $D$ and each of the $C_i$ are finite sets
in $\R$. By Lemma \ref{one}, C must be a finite set. Consider the cover
\[\{\pi_L(i) \: D \times_C C_i \rightarrow D | \, i \in I\},\] where
$\pi_L(i)$ is the obvious map and $g \circ \pi_L(i)$ factors
through $f_i$ via the canonical map $\pi_R(i)$. Now choose any $d \in
D$ and let $g(d) = c \in C$. Then there exists some $i$ such that
$f_i(c_i)=c$ for $c_i \in C_i$. Thus, $(d,c_i) \in D \times_C C_i$, so
that $\coprod_I D \times_C C_i \rightarrow D$ maps $(d,c_i)$ to $d$
and is therefore an epimorphism. This shows that $\{\pi_L(i) \}$ is in
$K(D)$.

{\it Case $\mathrm{(}$2$\, \mathrm{)}$}: Suppose that $D=\G$ and that each $C_i$ is a finite
set in $\R$. By Lemma \ref{one}, $C$ is a finite set and we identify it with
$\coprod_{i=1}^n\G/U_i$, where $U_i = {\G}_{x_i}$, the stabilizer of
$x_i$ in $\G$. The map $g$ is determined by $g(1) = \delta U_k$ for
some $\delta \in \G$ and some stabilizer $U_k$. Since $\coprod_I C_i
\rightarrow C$ is onto and $\mathrm{im}(g) = G/{U_k}$,
there exists some $c_l \in C_l$ such that $f_l(c_l) = U_k$. Since $C_l$
is a finite set, we can identify $c_l$ with some $\mu {\G}_z$, where $\mu
\in \G$ and ${\G}_z$ is the stabilizer of some element $z \in C_l$.

Then define the cover to be $\{\lambda \: \G \rightarrow
\G \}$, where $\lambda(\gamma) = \gamma \delta^{-1}$. Define $\alpha_l
\: \G  \rightarrow C_l$ to be the $\G$-equivariant map given by $1
\mapsto \mu {\G}_z$. By Lemma \ref{five}, $\alpha_l$ is continuous and is a
morphism in $\R$. Since $\lambda$ is a homeomorphism, the cover $\{\lambda\}$ is in $K(D)$. Now, \[(g \circ \lambda) (1) =
g(\delta^{-1}) = \delta^{-1} \cdot g(1) = U_k = \mu \cdot f_l( {\G}_z ) =
\mu \cdot f_l( \mu^{-1} \cdot \alpha_l(1) ) = (f_l \circ \alpha_l) (1).\] This shows
that $g \circ \lambda$ factors through $f_l$ via $\alpha_l$.

{\it Case $\mathrm{(}$3$\, \mathrm{)}$}: Suppose not all the $C_i$ are finite sets and that
$D=\G$. Also, assume that $C=\G$. This implies that $C_i=\G$ for all
$i \in I$. Choose any $k \in I$, let $\alpha_k = \mathrm{id}_{\G}$, and
define $\lambda \: \G \rightarrow \G$ to be multiplication on the
right by $f_k(1)g(1)^{-1}$. Then the diagram \[\xymatrix{ \G
\rto^{\mathrm{id}_{\G}} \dto_{\lambda} & \G \dto^{f_k} \\ \G \rto^{g} & \G}\]
is commutative, since \[(g \circ \lambda) (1) = g(f_k(1)g(1)^{-1})=
f_k(1)g(1)^{-1} \cdot g(1) = f_k(1) = (f_k \circ \alpha_k) (1).\] Thus,
$g \circ \lambda$ factors through $f_k$ via $\alpha_k$, so that the
stability axiom is verified by letting the covering family be $\{
\lambda \}.$

{\it Case $\mathrm{(}$4$\, \mathrm{)}$}: Suppose that not all the $C_i$ are finite sets, $D=\G$, and
$C$ is a finite set. With $C$ as in Lemma \ref{four}, let $g(1)=\delta U_k
\in C$, as in Case (2). Then there exists some $l$ such that $f_l(c_l) =
U_k$, for some $c_l \in C_l$. Now we consider two subcases. 
\par
{\it Case $\mathrm{(}$4a$\, \mathrm{)}$}: Suppose that $C_l$ is a
finite set in $\R$. Just as in Case (2), we construct maps $\lambda$
and $\alpha_l$, so that $g\circ \lambda$ factors through $f_l$ via
$\alpha_l$ and $\{\lambda \} \in K(D)$.
\par 
{\it Case $\mathrm{(}$4b$\, \mathrm{)}$}: Suppose that
$C_l = \G$. By $\G$-equivariance, $f_l(1)= c_l^{-1} U_k$. Then define $\lambda \: \G \rightarrow \G$
by $1 \mapsto \delta^{-1}$ and $\alpha_l \: \G \rightarrow \G$ by
$1 \mapsto c_l$. Then $g \circ \lambda$ factors through $f_l$
via $\alpha_l$, since \[(g \circ \lambda) (1) = g(\delta^{-1})=
\delta^{-1} \cdot g(1) = U_k =
f_l(c_l) = (f_l \circ \alpha_l)(1).\] Thus, the cover $\{
\lambda \}$, as a homeomorphism, is in $K(D)$. This completes Case
(4).

Now we consider the final possibility, {\it Case $\mathrm{(}$5$\, \mathrm{)}$}: suppose that
not all of the $C_i$ are finite sets and suppose that $D$ is a finite
set. This implies that $C$ is a finite set. This case is more
difficult than the others because the cover consists of more than one
morphism and it combines the previous constructions. For each $d \in
D$, we make a choice of some ${c_l} \in C_l$ for some $l$, such that
${c_l}$ is in the preimage of $g(d)$ under $\coprod C_i
\rightarrow C$. Then write $D = D_{df} \coprod D_G$, where
$D_{df}$ is the set of all $d$ such that the corresponding $C_l$
is in a finite set, and $D_G$ is the set of all $d$ such that
the corresponding $C_l = \G$. Now consider the cover
$\{h_d \: D_d \rightarrow D | \, d \in D =D_{df} \coprod
D_G\}$, where
\[
D_d = \begin{cases}
	D \times_C C_d & \text{if $d \in D_{df}$},\\
	\G & \text{if $d \in D_G$}.
\end{cases}
\] 
If $d \in D_{df}$, then $h_d = \pi_L$ and $\alpha_d \: D \times_C C_d
\rightarrow C_d$ is the canonical map $\pi_R$; it is clear that
the required square commutes. Now suppose $d \in D_G$. Then
there exists $c_l \in C_l = \G$ for some $l$, such that $g(d) =
f_l(c_l)$. We write $f_l(1) = \theta U_k \in C$ for
some $\theta \in \G$ and for some stabilizer $U_k$. Then we define
$\alpha_d \: \G \rightarrow C_l=\G$ by $1 \mapsto
\theta^{-1}$. Also, we define $h_d \: \G \rightarrow D$ by $1
\mapsto (\theta^{-1}c_l^{-1}) \cdot d$. Lemma \ref{five} shows that $h_d$ is
a morphism in $\R$. Then we have the required commutative diagram
\[
\xymatrix{ \G \rto^{\alpha_d} \dto_{h_d} & \G \dto^{f_l} 
\\ 
D
\rto^{g} & C,}
\] 
since
\begin{align*}
(g \circ h_d) (1) & = g((\theta^{-1}c_l^{-1}) \cdot d) =
(\theta^{-1}c_l^{-1}) \cdot g(d)\\ & = (\theta^{-1}c_l^{-1}) \cdot
f_l(c_l) = f_l(\theta^{-1}) =
(f_l \circ \alpha_d)(1).
\end{align*}

The only remaining detail is to show that $\{h_d\} \in K(D)$; that is,
we must show that $\phi \: \coprod_D D_d
\rightarrow D$ is an epimorphism. Let $d$ be any element in
$D$. Suppose $d \in D_{df}$. Then, using our choice above, there
exists some $c_l \in C_l$, a finite set for some $l$, such that
$f_l(c_l)=g(d)$. Then $(d, c_l) \in D \times_C C_l$
and $\phi(d, c_l) = \pi_L(d, c_l) = d$. Now suppose $d \in
D_G$. With $c_l$ and $\theta$ as above, $c_l \theta \in D_d = \G$ and
$\phi(c_l\theta)= h_d(c_l\theta) = (c_l\theta) \cdot h_d(1)
= d$. Therefore, $\phi$ is an epimorphism.
\epf
\section{The site $\R$ does not have the canonical topology}
\par
Now that we have established that $\R$ is a site with pretopology
$\mathcal{K}$, we begin working to show that, contrary to what typically
happens with this pretopology, it does not give the canonical
topology. We start with a definition.
\begin{definition}
If $T$ is some collection of morphisms with codomain $C$, where $C$ is
an object in the category $\mathcal{C}$, then $(T)$ denotes the sieve
generated by $T$. Thus, \[(T) = \{ f \circ g | \, f \in T, \ \mathrm{dom}(f) = 
\mathrm{cod}(g) \}.\]
\end{definition}
\lemma Let $K$ be a pretopology on a category $\mathcal{C}$. Let
$J$ be the Grothendieck topology generated by $K$. Then for any $C \in
\mathcal{C}$, $J(C)$ consists exactly of all $(R) \cup (T)$ such that
$R \in K(C)$ and $T$ is some collection of morphisms with codomain
$C$.
\endlemma
\pf
Let $S$ be a covering sieve of $C$. Then there exists some $R \in
K(C)$ such that $R \subset S$. We will prove that $S = (R) \cup (S)$,
verifying the forward inclusion. To prove equality it suffices to
show that $(R) \cup (S) \subset S$. If $f \in (R)$, $f = g \circ h$
for some $g \in R$ and some $h$ with $\mathrm{dom}(g) = \mathrm{cod}(h)$. 
Since $g \in S$,
$f \in S$. Similarly, if $f \in (S)$, then $f \in S$. Now consider any family of morphisms $(R) \cup (T)$
as described in the statement of the lemma. Since $R \subset (R) \cup
(T)$, $(R) \cup (T) \in J(C)$ if it is a sieve. Since $(R)$ and $(T)$
are sieves, it is clear that $(R) \cup (T)$ is also a sieve.
\epf
\par
This result is useful for understanding the topology of a site, when
the site is defined in terms of a pretopology. For example,
$G \negthinspace - \negthinspace \mathbf{Sets}_{df}$ is a site by the pretopology 
$\mathcal{K}$ and its
category of sheaves of sets is equivalent to the category of sheaves
on the site $S(G)$ consisting of quotients of $G$ by open subgroups
(the morphisms are the $G$-equivariant maps), where $S(G)$ is given
the atomic topology (see \cite[Chapter 3, Section 9]{MacLane}). Thus, one might ask if $G \negthinspace - \negthinspace \mathbf{Sets}_{df}$ also has the
atomic topology. However, the lemma allows us to see that $\mathcal{K}$
generates a topology that is coarser than the atomic topology. To see this, let $X = G/U$ and $Y = G/U \coprod G/U$, where $U$
is a proper open subgroup of $G$. (Since $G$ is an infinite profinite
group, the canonical way of writing $G$ as an inverse limit
guarantees the existence of such a $U$.) We define $f \: X
\rightarrow Y$ by $f(U) = U$, where $U$ lives in the factor
on the left; $f$ is the left inclusion. Now consider the sieve $S =
(\{f\})$. Clearly, $S$ does not contain an epimorphic cover, since
$\mathrm{im}(\coprod_{g \in S} (\mathrm{dom}(g)) \rightarrow Y) = G/U$. The lemma
indicates that every sieve of $G \negthinspace - \negthinspace \mathbf{Sets}_{df}$ must contain an epimorphic
cover, so that $S$ is not a sieve for $Y$ in the topology generated by
$\mathcal{K}$. 
\par
Now we consider the site $R^+_G$ with the pretopology $\mathcal{K}$ of
epimorphic covers. We use $\mathrm{Hom}_G(X,Y)$ to denote continuous
$G$-equivariant maps between continuous $G$-sets $X$ and $Y$. Recall that a presheaf of sets P on a site $(\mathcal{C}, J)$ is a sheaf,
if for each object $C \in \mathcal{C}$ and each covering sieve $S \in
J(C)$, the diagram
\[\xymatrix{
P(C)
\ar[r]^-{e} 
&
\prod_{f \in S} P(\mathrm{dom}(f)) 
\ar@<1ex>[r]^-{p} \ar@<-1ex>[r]_-{a} &
\prod
P(\mathrm{dom}(g))}
\]
is an equalizer of sets, where the second product is over all $f, g,$ with 
$f \in S,$ $\mathrm{dom}(f) = \mathrm{cod}(g)$. Here, $e$ is the map $e(x) = \{
P(f)(x) \}_f$, $p$ is given by \[\{x^f\}_f \mapsto \{x^{fg}\}_{f,g},\]
and $a$ is given by \[\{x^f\}_f \mapsto \{P(g)(x^{f})\}_{f,g} =
\{x^{f} \circ g \}_{f,g}.\] 
\par
Recall that a representable presheaf of $\R$ is any presheaf which, up
to isomorphism, has the form of $\mathrm{Hom}_G(-,C)$ for some $C \in
\R$. Also, the Yoneda embedding \[\R \rightarrow \mathbf{Sets}^{({\R})^\mathrm{op}},
\ \ \ C \mapsto \mathrm{Hom}_G(-,C)\] is a full and faithful functor, so
that one can identify $C$ with an object of $\mathbf{Sets}^{({\R})^\mathrm{op}}$. We
now consider which objects of $\R$ yield sheaves of sets on $\R$.

\par
Noting that the empty set is a discrete $G$-set, we have

\lemma \label{stupidsheaf}
The presheaf $\mathrm{Hom}_G(-,\varnothing)$ is a sheaf of sets on the site
$\R$.
\endlemma

\pf
Let $\bullet: \varnothing \rightarrow X$ denote the vacuous
map, for any $X \in \R$. Since $\bullet: \varnothing \rightarrow \varnothing$ is
vacuously an epimorphism, $\{\bullet\}$ is the unique covering sieve
for $\varnothing$. Let $C= \varnothing$. Then the desired equalizer
diagram has the form 
\[
\xymatrix{\mathrm{Hom}_G(\varnothing, \varnothing) =
 \{\bullet\}
\ar[r]^-{e} 
&
 \{\bullet\}
\ar@<1ex>[r]^-{p} \ar@<-1ex>[r]_-{a} &
\{\bullet\}. \\}\] It is clear that this is an equalizer diagram. 
\par
Now
let $C$ be a nonempty finite set in $G\negthinspace-\negthinspace \mathbf{Sets}_{df}$. Let $S$ be any
covering sieve of $C$. There must exist a morphism in $S$ with domain
equal to a nonempty object in $\R$. Therefore, since $\varnothing \times Z = \varnothing$ for any space $Z$, we have
\[
\xymatrix{\mathrm{Hom}_G(C, \varnothing) = \varnothing
\ar[r]^-{e} 
&
 \varnothing
\ar@<1ex>[r]^-{p} \ar@<-1ex>[r]_-{a} &
\varnothing. \\}\]
Since the equalizer must exist and the vacuous map
$\bullet \spc : \varnothing \rightarrow \varnothing$ is the unique
map with codomain $\varnothing$, this must be an equalizer
diagram.
\par
Finally, letting $C = G$, we get
\[
\xymatrix{
\mathrm{Hom}_G(G, \varnothing) = \varnothing
\ar[r]^-{e} 
&
\prod_{f \in \mathrm{Hom}_G(G,G)} \varnothing = \varnothing
\ar@<1ex>[r]^-{p} \ar@<-1ex>[r]_-{a} &
\varnothing. \\}\] Again, this is an equalizer diagram. 
\epf
\par
To prove the next theorem, we need the following lemma.
\lemma
If $G$ is a compact topological group, $U$ an open subgroup of $G$, and $X \neq
\varnothing$ a finite discrete $G$-set, then \[\mathrm{Hom}_G(G/U,X) \cong \{ x
\in X | \, U < G_x \},\] where $G_x$ is the stabilizer of $x$ in $G$.
\endlemma

\pf
Let $f \spc : G/U \rightarrow X$. It is clear that $f$ is
$G$-equivariant if and only if it is completely
determined by $f(U)$ in the obvious way. Since $U$ is an open subgroup, it has finite
index in $G$, so that $G/U$ is a discrete space. Thus, any
$G$-equivariant map $G/U \rightarrow X$ is continuous. The key is
that $f$ is well-defined if and only if $U < G_{f(U)}$. To see this, first
assume that $f$ is well-defined; let $\gamma \in U$. Then $\gamma U =
U$, so that $\gamma \cdot f(U) = f(\gamma U) =
f(U)$. Hence, $\gamma \in G_{f(U)}$ and $U < G_{f(U)}$. Now suppose
that $U < G_{f(U)}$ and take any $\gamma U = \delta U$. This implies
that $\gamma^{-1}\delta \in U$ and hence, in $G_{f(U)}$. Thus,
$(\gamma^{-1}\delta) \cdot f(U) = f(U)$, so that $\gamma \cdot f(U) =
\delta \cdot f(U)$. Equivariance gives $f(\gamma U) = f(\delta U)$ and
$f$ is well-defined. Thus, \[\mathrm{Hom}_G(G/U,X) \cong \{ f(U) \in X | \, U <
G_{f(U)} \}.\]
\epf 

\par
Henceforth, let $\mathcal{J}$ denote the topology of $\R$ generated by $\mathcal{K}$.

\theorem \label{notsub}
Let $X$ be any object in $\R$ that is not a finite 
discrete trivial $G$-set, where $G$
is an infinite profinite group. Then the
presheaf $\mathrm{Hom}_G( -, X)$ is not a sheaf of sets on the site $R^+_G$.
\endtheorem

\pf
Suppose $\mathrm{Hom}_G(-,X)$ is a sheaf of sets on the site $\R$. The
equalizer condition says that for every object $C \in \R$ and for
every covering sieve $S \in \mathcal{J}(C)$, \[\mathrm{Hom}_G(C,X) \cong \{ \{h^f\}_f |
\, h^{fg} = h^f \circ g, f, g, f \in S, \mathrm{dom}(f) = \mathrm{cod}(g) \},\] where for
$f \in S$, $h^f \in \mathrm{Hom}_G(\mathrm{dom}(f), X)$. 
We will construct an example of
some $C$ and $S$ such that this sheaf condition
fails to be true with $X$ as above. 
\par
Let $C \in G \negthinspace - \negthinspace \mathbf{Sets}_{df}$; we identify $C$ with $\coprod^n_{i=1}
G/{U_i}$, where each $U_i$ is an open subgroup of $G$. For each $i$,
define $f_i \: G \rightarrow C$ by $1 \mapsto U_i$. Thus, $\mathrm{im}(f_i) =
G/{U_i}$ and $\{f_i\}$ is an epimorphic cover of $C$. The preceding lemma tells
us that $S = (\{f_i\})$ is a covering sieve of $C$. For this $S$, we
will examine the sheaf condition. Let $S = S' \cup S''$, where $S'
=\{f_i\}$ and $S''$ is the complement of $S'$ in $S$. Thus, every $k
\in S''$ has the form $k = f_i \circ g$ for some $g$ with $\mathrm{dom}(f_i) =
\mathrm{cod}(g)$. Then 
\begin{align*}
\{ \{h^f\}_f |
\, h^{fg} & = h^f \circ g, f, g, f \in S, \mathrm{dom}(f) = \mathrm{cod}(g) \} \\
& = \{ \{h^{f_i}\}_{f_i} \times \{h^k\}_{k \in S''} |
\, h^{fg} = h^f \circ g, f, g, f \in S, \mathrm{dom}(f) = \mathrm{cod}(g) \}\\ & =
\{ \{h^{f_i}\}_{f_i} \times \{h^{f_i} \circ g\}_{f_i \circ g \in S''} |
\, h^{fg} = h^f \circ g, f, g, f \in S, \mathrm{dom}(f) = \mathrm{cod}(g) \} \\ & = 
\{ \{h^{f_i}\}_{f_i} \times \{h^{f_i} \circ g\}_{f_i \circ g \in S''} |
\, h^{f_i} \negthinspace \in \negthinspace \mathrm{Hom}_G(G,X), f_i \negthinspace \in \negthinspace S', g, \mathrm{dom}(f_i) = \mathrm{cod}(g) \}.
\end{align*}
We verify the last equality. Suppose $h^{f_i}$ is any morphism in
$\mathrm{Hom}_G(G,X)$. Now take any $f$ and $g$ with $f \in S$ and 
$\mathrm{dom}(f) = \mathrm{cod}(g)$. If $f = f_i \in S'$, then $h^f \circ g = h^{f_i} \circ g = h^{f_i
\circ g} = h^{f \circ g},$ by construction. Now suppose $f \in
S''$. Then $f = f_i \circ k$ for some $k \: G \rightarrow G$. Thus,
\[h^{fg} = h^{f_i \circ (k \circ g)} = (h^{f_i} \circ k) \circ g =
h^{f_i \circ k} \circ g = h^{f} \circ g.\] Since $h^{f_i} \circ g$ is
determined by $h^{f_i}$ and $f_i \circ g$, we see that the set 
\[
\{ \{h^{f_i}\}_{f_i} \times \{h^{f_i} \circ g\}_{f_i \circ g \in S''} |
\, h^{f_i} \in \mathrm{Hom}_G(G,X), f_i \in S', g, \mathrm{dom}(f_i) = \mathrm{cod}(g) \}
\] 
is isomorphic to the set 
\[
\{ \{h^{f_i}\}_{f_i} |
\, h^{f_i} \in \mathrm{Hom}_G(G,X), f_i \in S' \} = \mathrm{Hom}_G(G,X)^n,
\] where
$\mathrm{Hom}_G(G,X)^n$ is the $n$-fold Cartesian product of $\mathrm{Hom}_G(G,X)$. Now, there is an isomorphism 
$\mathrm{Hom}_G(G,X)^n \cong X^n.$ Therefore, for $\mathrm{Hom}_G(-,X)$ to be a sheaf, it must be that 
$\mathrm{Hom}_G(C,X)
\cong X^n$ for every $C \in G-\mathbf{Sets}_{df}.$ If $X=G$ and $C \neq
\varnothing$ is in $G-\mathbf{Sets}_{df}$, then $\mathrm{Hom}_G(C,G) =
\varnothing$, whereas, since $|C| \geq 1$, $n \geq 1$ and $X^n =
G^n$. Thus, $\mathrm{Hom}_G(-,G)$ is not a
sheaf. 
\par
Now we consider $X \neq G$ and assume that $\mathrm{Hom}_G(C,X)
\cong X^n$ for every $C \in G\negthinspace-\negthinspace \mathbf{Sets}_{df}.$ This implies that 
\begin{align*}
X^n & \cong
\mathrm{Hom}_G(C,X) \cong \mathrm{Hom}_G(\textstyle{\coprod}^n_{i=1} 
G/{U_i},X) \\ & 
\cong \textstyle{\prod}^n_{i=1} \mathrm{Hom}_G(G/{U_i},X) \cong 
\textstyle{\prod}^n_{i=1} \{ x
\in X | \, U_i < G_x \} \subset X^n.\end{align*}
Therefore, it must be that $\{ x
\in X | \, U_i < G_x \} = X$, for all $i = 1, ..., n$. Thus, $U_i <
G_x$ for all $x \in X$ and each $i$. Now let us write $X \cong
\coprod^m_{j=1} G/{G_{x_j}}$, where each $x_j$ is a representative
from a distinct orbit of $X$. Let $C$ be a trivial $G$-set so
that every stabilizer of $c \in C$ in $G$ is equal to $G$. This implies that
$G < G_{x_j}$ for all $j$. Thus, each $G_{x_j} = G$. This indicates
that $X$ must be a trivial $G$-set. This contradiction shows that
every $X$ violates the sheaf condition for some $C$ and $S$.
\epf 
\par 
This result immediately yields
\corollary
For an infinite profinite group $G$, the site $\R$ with the pretopology
$\mathcal{K}$ of epimorphic covers is not subcanonical.
\endcorollary
\pf
There exists a proper open subgroup $U$ of $G$ satisfying $[G:U] >
1$. Thus, the representable presheaves $\mathrm{Hom}_G(-,G)$ and $\mathrm{Hom}_G(-,
\coprod^n_{i=1} G/U)$, for any $n \geq 1$, are not sheaves.
\epf
Since a canonical topology is, by definition, subcanonical, we obtain
\corollary
For an infinite profinite group $G$, the site $\R$, with the pretopology $\mathcal{K}$, 
is not canonical.
\endcorollary
\par
The next result is an elementary fact about profinite groups that
helps us understand ``how often'' representable presheaves fail to be
sheaves in $\R$ and what such ``failing'' presheaves can look
like, based on what we know from Theorem \ref{notsub}.
\lemma
If $G$ is an infinite profinite group, then $G$ contains an infinite
number of distinct proper open subgroups.
\endlemma
\pf
We have already seen that $G$ has at least one proper open
subgroup. Suppose that $G$ has only a finite number of distinct proper
open subgroups. Then $G$ has a finite number of distinct proper open
normal subgroups $N_1, ..., N_k$. Since $G$ is profinite,
$N = \bigcap_{i=1}^k N_i = \{1\}$. Because $N$ is an open subgroup with
finite index, it has uncountable order. This contradiction gives the conclusion.
\epf
\begin{remark}
Since any topology finer than $\mathcal{J}$ would contain the covering
sieve $(\{f_i\})$ that was the key to Theorem \ref{notsub}, no topology finer
than $\mathcal{J}$ can be subcanonical.
\end{remark}
\section{Presheaves of spectra on the site $\R$}\label{presheaves}
\par
Let $\mathbf{Ab}$ be the category of abelian groups, and let $\mathrm{Spt}$ denote the model category of Bousfield-Friedlander spectra of pointed simplicial sets. We refer to the objects of $\mathrm{Spt}$ as simply ``spectra." Now that $\R$ is a site, we can consider the category $\mathrm{PreSpt}(\R)$ of 
presheaves of spectra on the site $\R$. By applying the work of Jardine (\cite{Jardinecanada}, \cite[Section 2.3]{Jardine}), $\mathrm{PreSpt}(\R)$ is a model category. We recall the critical definitions that give the model category structure and then we state Jardine's result, when it is applied to $R^+_G$. 
\begin{definition} 
Let $P \: (R^+_G)^\mathrm{op} \rightarrow \mathrm{Spt}$ be a presheaf of spectra. Then, for 
each $n \in \mathbb{Z}$,
\[\pi_n(P) \: (R^+_G)^\mathrm{op} \rightarrow \mathbf{Ab}, \ \ \ C \mapsto \pi_n(P(C)),\] is a presheaf of abelian groups. Then the associated sheaf $\tilde{\pi}_n(P)$ of abelian groups is the 
sheafification of $\pi_n(P).$
\end{definition}
\par
Let $f \: P \rightarrow Q$ be a morphism of presheaves of spectra on $R^+_G$. Then $f$ is 
a {\it weak equivalence} if the induced map $\tilde{\pi}_n(P) \rightarrow \tilde{\pi}_n(Q)$ of sheaves 
is an isomorphism, for all $n \in \mathbb{Z}$. The map $f$ is a {\it cofibration} if $f(C)$ 
is a cofibration of spectra, for all $C \in \R$. Also, $f$ is a {\it global fibration} if 
$f$ has the right lifting property with respect to all morphisms which are weak equivalences and cofibrations. 
\theorem[{\cite[Theorem 2.34]{Jardine}}]
The category $\mathrm{PreSpt}(\R)$, together with the classes of weak equivalences, cofibrations, and global fibrations, is a model category.
\endtheorem
\par
Now we give some interesting examples of presheaves of spectra on the site $R^+_G$.
\begin{example}
In the Introduction, we saw that the Devinatz-Hopkins functor $\mathbf{F}$ is an example of an object in $\mathrm{PreSpt}(R^+_{G_n})$.  
\end{example}
\par
For the next example, if $X$ is a spectrum, then, for each $k \geq 0$, we let $X_k$ be the $k$th pointed simplicial set constituting $X$, and, for each $l \geq 0$, $X_{k,l}$ is the pointed set of $l$-simplices of $X_k$. 
\begin{example} 
Let $X$ be a discrete $G$-spectrum (see \cite{cts} for a definition of this term), so that 
each $X_{k,l}$ is a pointed discrete $G$-set. If $C \in R^+_G$, then let $\mathrm{Hom}_G(C,X)$ be the spectrum, such that 
\[\mathrm{Hom}_G(C,X)_k = \mathrm{Hom}_G(C,X_k),\] where 
\[\mathrm{Hom}_G(C,X)_{k,l} = \mathrm{Hom}_G(C,X_k)_l = \mathrm{Hom}_G(C,X_{k,l}).\] Above, 
the set $X_{k,l}$ is given the discrete topology, since it is naturally a discrete $G$-set. Then 
$\mathrm{Hom}_G(-,X)$ is an object in $\mathrm{PreSpt}(R^+_G)$. It is easy to see that 
if $U$ is an open subgroup of $G$, then $\mathrm{Hom}_G(G/U,X) \cong X^U$, the $U$-fixed point spectrum of $X$. Also, $\mathrm{Hom}_G(G,X) \cong X$.
\end{example}
\par
Now we recall part of \cite[Proposition 3.3.1]{joint}, since this result (and its corollary) will be helpful in our next example. We note that this result is only a slight extension of \cite[Remark 6.26]{Jardine}: if $U$ is normal in $G$, then the lemma below is an immediate consequence of Jardine's remark.  
\lemma
Let $X$ be a discrete $G$-spectrum. Also, let $f \: X \rightarrow X_{f,G}$ be a trivial cofibration, such that $X_{f,G}$ is fibrant, where all this takes place in the model 
category of discrete $G$-spectra $\mathrm{(}$see $\mathrm{\cite{cts})}$. If $U$ is an open subgroup of 
$G$, then $X_{f,G}$ is fibrant in the model category of discrete $U$-spectra. 
\endlemma
\corollary
Let $X$ and $U$ be as in the preceding lemma. Then $X^{hU} = (X_{f,G})^U$.
\endcorollary
\pf
Let $f$ be as in the above lemma. Since $f$ is $G$-equivariant, it is $U$-equivariant. Also, since $f$ is a trivial cofibration in the model category of discrete $G$-spectra, it is a trivial cofibration in the model category of spectra. The preceding two facts imply that 
$f$ is a trivial cofibration in the model category of discrete $U$-spectra. By the lemma, $X_{f,G}$ is fibrant in this model category. Thus, $X^{hU} = (X_{f,G})^U$.
\epf
\begin{example}\label{exone}
Let $X$ be a discrete $G$-spectrum. Then $\mathrm{Hom}_G(-,X_{f,G})$ is a presheaf in 
$\mathrm{PreSpt}(\R)$. In particular, notice that \[\mathrm{Hom}_G(G/U,X_{f,G}) \cong 
(X_{f,G})^U = X^{hU}\] and \[\mathrm{Hom}_G(G,X_{f,G}) \cong X_{f,G} \simeq X.\]
\end{example}
\begin{example}\label{extwo}
For any unfamiliar concepts in this example, we refer the reader to \cite{cts}. Let $Z = 
\holim_i Z_i$ be a continuous $G$-spectrum, so that $\{Z_i\}_{i \geq 0}$ is a tower of discrete $G$-spectra, such that each $Z_i$ is a fibrant spectrum. Then 
\[P(-) = \holim_i \mathrm{Hom}_G(-,(Z_i)_{f,G}) \in \mathrm{PreSpt}(\R),\] where 
\[P(G/U) \cong \holim_i ((Z_i)_{f,G})^U = \holim_i (Z_i)^{hU} = Z^{hU}\] and 
\[P(G) \cong \holim_i (Z_i)_{f,G} \simeq Z.\]
\end{example} 


\refs

\bibitem
[Behrens and Davis, 2005]{joint}
Mark~J. Behrens and Daniel~G. Davis, \emph{The homotopy fixed point spectra of
  profinite {G}alois extensions}, 33 pp., current version of manuscript is
  available at http://dgdavis.web.wesleyan.edu, 2005.

\bibitem
[Davis, 2006]{cts}
Daniel~G. Davis, \emph{Homotopy fixed points for {$L_{K(n)}(E_n \wedge X)$}
  using the continuous action}, J. Pure Appl. Algebra \textbf{206} (2006),
  no.~3, 322--354.

\bibitem
[Demazure, 1970]{Demazure}
Michel Demazure, \emph{Topologies et faisceaux}, Sch\'emas en groupes.
  {S}{G}{A} {3}, {I}: {P}ropri\'et\'es g\'en\'erales des sch\'emas en groupes,
  Springer-Verlag, Berlin, 1970, pp.~xv+564.

\bibitem
[Devinatz and Hopkins, 1995]{DevinatzHopkins}
Ethan~S. Devinatz and Michael~J. Hopkins, \emph{The action of the {M}orava
  stabilizer group on the {L}ubin-{T}ate moduli space of lifts}, Amer. J. Math.
  \textbf{117} (1995), no.~3, 669--710.

\bibitem
[Devinatz and Hopkins, 2004]{DH}
Ethan~S. Devinatz and Michael~J. Hopkins, \emph{Homotopy fixed point spectra for closed subgroups of the
  {M}orava stabilizer groups}, Topology \textbf{43} (2004), no.~1, 1--47.

\bibitem
[Elmendorf et. al., 1997]{EKMM}
A.~D. Elmendorf, I.~Kriz, M.~A. Mandell, and J.~P. May, \emph{Rings, modules,
  and algebras in stable homotopy theory}, American Mathematical Society,
  Providence, RI, 1997, With an appendix by M. Cole.

\bibitem
[Jardine, 1987]{Jardinecanada}
J.~F. Jardine, \emph{Stable homotopy theory of simplicial presheaves}, Canad.
  J. Math. \textbf{39} (1987), no.~3, 733--747.

\bibitem
[Jardine, 1997]{Jardine}J.~F. Jardine, \emph{Generalized \'etale cohomology theories}, Birkh\"auser Verlag,
  Basel, 1997.

\bibitem
[Mac~Lane and Moerdijk, 1994]{MacLane}
Saunders Mac~Lane and Ieke Moerdijk, \emph{Sheaves in geometry and logic},
  Springer-Verlag, New York, 1994, A first introduction to topos theory,
  Corrected reprint of the 1992 edition.

\bibitem
[Rudyak, 1998]{Rudyak}
Yuli~B. Rudyak, \emph{On {T}hom spectra, orientability, and cobordism},
  Springer Monographs in Mathematics, Springer-Verlag, Berlin, 1998, With a
  foreword by Haynes Miller.
  
\endrefs

\end{document}